\documentclass[envcountsame]{llncs}

\usepackage{amssymb}
\usepackage{a4}
\usepackage{QED}

\newcommand {\es}{\emptyset}
\newcommand {\ra}{\tr}
\newcommand {\ras}{\tr^*}

\newcommand {\G}{\Gamma}
\newcommand {\m}{\mu}

\newcommand {\be}{\beta}
\newcommand {\al}{\alpha}
\newcommand {\sig}{\sigma}
\newcommand {\la}{\lambda}
\newcommand {\ptv}{\; ..., \;}
\newcommand {\pts}{ \; ... \; }
\newcommand {\peq}{\preceq}

\def\b{\beta}

\def\G{\Gamma}

\def\l{\lambda}

\def\m{\mu}

\def\f{\rightarrow}
\def\tr{\triangleright}

\def\v{\vdash}

\def\<{\langle}
\def\>{\rangle}
\def\F{\displaystyle\frac}

\begin{document}

\begin{center}

{\Large\bf
Arithmetical proofs of strong normalization results for the symmetric $\lambda \mu$-calculus}\\[1cm]
\end{center}

\begin{center}
{\bf Ren\'e David and Karim Nour}\\
Laboratoire de Mathé\'ematiques\\
 Universit\'e de Savoie\\
73376 Le Bourget du Lac. France\\
e-mail: \{david,nour\}@univ-savoie.fr\\[1cm]
\end{center}

\begin{abstract}
The  symmetric $\lambda \mu$-calculus is the  $\lambda
\mu$-calculus introduced by Parigot in which the reduction rule
$\m'$, which is the symmetric of $\m$, is added. We give {\em
arithmetical} proofs of some strong normalization results for this
calculus. We show (this is a new result) that the
$\m\m'$-reduction is strongly normalizing for the un-typed
calculus. We also show the strong normalization of the
$\be\m\m'$-reduction for the typed calculus: this was already
known but the previous proofs use candidates of reducibility where
the interpretation of a type was defined as the fix point of some
increasing operator and thus, were highly non arithmetical.
\end{abstract}

\section{Introduction}

Since it has been understood that  the Curry-Howard isomorphism
relating proofs and programs  can be extended to classical logic,
various  systems have been introduced:  the $\l_c$-calculus
(Krivine \cite{Kri}), the $\la_{exn}$-calculus (de Groote
\cite{deG4}), the $\l \mu$-calculus (Parigot \cite{Par1}), the
$\lambda^{Sym}$-calculus (Barbanera \& Berardi \cite{BaBe}), the
$\lambda_{\Delta}$-calculus (Rehof \& Sorensen \cite{ReSo}), the
$\overline{\lambda}\mu\tilde{\mu}$-calculus (Curien \& Herbelin
\cite{cuhe}),  ...

The first calculus which respects the intrinsic symmetry of
classical logic is $\lambda^{Sym}$. It is somehow different from
the previous calculi since the main connector is not the arrow as
usual but the connectors {\em or} and {\em and}. The symmetry of
the calculus comes from the de Morgan laws.

The second calculus respecting this symmetry has been
$\overline{\lambda}\mu\tilde{\mu}$. The logical part is the
(classical) sequent calculus instead of natural deduction.

Natural deduction is not, intrinsically, symmetric but Parigot has
introduced the so called {\em Free deduction} \cite{Pa01} which is
completely symmetric. The  $\l \mu$-calculus  comes from there. To
get a confluent calculus he had, in his terminology, to fix the
inputs on the left. To keep the symmetry, it is enough to keep the
same terms and to add a new reduction rule (called the
$\mu'$-reduction) which is the symmetric rule of the
$\mu$-reduction and also corresponds  to the elimination of a
 cut. We get then a symmetric calculus that is called the
{\em symmetric $\lambda \mu$-calculus}.

The $\m'$-reduction has been considered by Parigot for the
following reasons. The $\l \mu$-calculus (with the
$\beta$-reduction and the $\mu$-reduction) has good properties :
confluence in the un-typed version, subject reduction and strong
normalization in the typed calculus. But this system has, from a
computer science point of view, a drawback: the unicity of the
representation of data is lost. It is known that, in the
$\la$-calculus, any term of type $N$ (the usual type for the
integers) is $\be$-equivalent to a Church integer. This no more
true in the $\la\m$-calculus and we can find normal terms of type
$N$ that are not Church integers. Parigot has remarked that by
adding the $\mu'$-reduction and some simplification rules the
unicity of the representation of data is recovered and subject
reduction is preserved, at least for the simply typed system, even
though the confluence is lost.

Barbanera \& Berardi proved the strong normalization of the
$\lambda^{Sym}$-calculus by using candidates of reducibility but,
unlike the usual construction (for example for Girard's system
$F$), the definition of the interpretation of a type needs a
rather complex fix-point operation. Yamagata \cite{Yam} has used
the same technic to prove the strong normalization of the symmetric
$\la\m$-calculus where the types are those of system $F$ and
Parigot, again using the same ideas, has extended Barbanera \&
Berardi's result to a logic with second order quantification.
These proofs are thus highly non arithmetical.

We consider here  the $\lambda \mu$-calculus with the rules $\be$,
$\m$ and $\m'$. It was known that, for the un-typed calculus, the
$\m$-reduction is strongly normalizing (see \cite{Py}) but the
strong normalization of the $\m\m'$-reduction for the un-typed
calculus was an open problem  raised long ago by Parigot. We give
here a proof of this result. Studying this reduction by itself is
interesting since a $\m$ (or $\m'$)-reduction can be seen as  a
way ``to put the arguments of the $\m$ where they are used'' and
it is useful to know that this is terminating. We also give an
{\em arithmetical} proof of the strong normalization of the
$\be\m\m'$-reduction for the simply typed calculus. We finally
show (this is also a new result) that, in the un-typed calculus,
if $M_1,...,M_n$ are strongly normalizing for the
$\beta\m\m'$-reduction, then so is $(x \; M_1 \pts M_n)$.

The proofs of strong normalization that are given here are
extensions of the ones given by the first author for the simply
typed $\l$-calculus. This proof can be found either in \cite{dav1}
(where it appears among many other things) or as a simple
unpublished note on the web page of the first author
 (\verb www.lama.univ-savoie.fr/~david ~).

The same proofs can be done for the
$\overline{\lambda}\mu\tilde{\mu}$-calculus and these proofs are,
in fact, much simpler for this calculus since some difficult
problems that appear in the $\l\m$-calculus do not appear in the
$\overline{\lambda}\mu\tilde{\mu}$-calculus: this is mainly due to
the fact that, in the latter, there is a right-hand side and a
left-hand side (the terms and the environments) whereas, in the
$\l\m$-calculus, this distinction is impossible since a term on
the right of an application can go on the left of an application
after some reductions. The proof of the strong normalization of
the $\mu\tilde{\mu}$-reduction can be found in \cite{polo}. The
proof is done (by using candidates of reducibility and a fix point
operator) for a typed calculus but, in fact, since the type system
is such that every term is typable, the result is valid for every
term. A proof of the strong normalization of the
$\overline{\lambda}\mu\tilde{\mu}$-typed calculus (again using
candidates of reducibility and a fix point operator) can also be
found there. Due to the lack of space, we do not give our proofs
of these results here but they will appear in \cite{craco1}.

The paper is organized as follows. In section \ref{s2} we give the
syntax of the terms and the reduction rules.  An arithmetical
proof of  strong normalization is given in section \ref{mm'} for
the $\m\m'$-reduction of the un-typed calculus and, in section
\ref{beta}, for the $\beta\m\m'$-reduction of the simply typed
calculus. In section \ref{s3}, we give an example showing that the
proofs of strong normalization using candidates of reducibility
{\em must} somehow be different from the usual ones and we show
that, in the un-typed calculus, if $M_1, ..., M_n$ are strongly
normalizing for the $\beta\m\m'$-reduction, then so is $(x \; M_1
...\; M_n)$. We conclude with some future work.

\section{The symmetric $\l\m$-calculus }\label{s2}

\subsection{The un-typed calculus}

The set (denoted as ${\cal T}$) of $\l\m$-terms or simply terms is
defined by the following grammar where $x,y,...$ are
$\lambda$-variables and $\al, \be, ...$ are $\mu$-variables:
$$
{\cal T} ::= x \mid \l x {\cal T} \mid ({\cal T} \; {\cal T}) \mid
\mu \al {\cal T} \mid   (\al \; {\cal T})
$$

Note that we adopt here a more liberal syntax (also called de
Groote's calculus) than in the original calculus since we do not
ask that a $\m \al$ is immediately followed by a $(\be \; M)$
(denoted $[\be] M$ in Parigot's notation).

\begin{definition}\rm
Let $M$ be a term.
\begin{enumerate}
\item $cxty(M)$ is the number of  symbols occurring in $M$.

\item
  We denote by  $N \leq M$ (resp. $N < M$) the fact that $N$ is a sub-term
  (resp. a strict sub-term) of $M$.
  \item If $\overrightarrow{P}$ is a sequence $P_1,...,P_n$ of terms,  $(M \; \overrightarrow{P})$
  will denote $(M \; P_1 \pts P_n)$.
\end{enumerate}
\end{definition}

\subsection{The typed calculus}

The types  are those of the simply typed $\l\m$-calculus i.e. are
built from atomic formulas and the constant symbol $\perp$ with
the connector $\rightarrow$. As usual $\neg A$ is an abbreviation
for $A \rightarrow \perp$.

 The typing rules are given by figure 1 below
where $\G$ is a context, i.e. a set of declarations of the form $x
: A$ and $\al : \neg A$ where $x$ is a $\l$ (or intuitionistic)
variable, $\al$ is a $\m$ (or classical) variable and $A$ is a
formula.

\begin{center}
$\F{}{\G , x : A \v x : A} \, ax$

\medskip
$\F{\G, x : A \v M : B} {\G \v \l x M : A \f B} \, \f_i$
\hspace{0.5cm}  $\F{\G \v M : A \f B \quad \G \v N : A} {\G \v (M
\; N): B }\, \f_e$

\medskip

$\F{\G , \al : \neg A  \v M : \bot} {\G \v \mu \al M : A }  \,
\bot_e$ \hspace{0.5cm} $\F{\G , \al : \neg A  \v M : A} {\G, \al :
\neg A  \v ( \al \; M) : \bot } \, \bot_i$
\medskip

Figure 1.
\end{center}

Note that, here, we also have changed Parigot's notation but these
typing rules are those of his classical natural deduction. Instead
of writing
$$M : (A_1^{x_1 }, ..., A_n^{x_n }  \vdash   B,  C_1^{\al_1},
..., C_m^{\al_m })$$ we have written
$$x_1 : A_1, ..., x_n : A_n,
\al_1: \neg C_1, ..., \al_m : \neg C_m \vdash M : B$$

\begin{definition}\rm
Let $A$ be a type. We  denote by $lg(A)$ the number of arrows in
$A$.
\end{definition}

\subsection{The reduction rules}

The cut-elimination procedure (on the logical side) corresponds to
the reduction rules (on the terms) given below. There are three
kinds of cuts.

\begin{itemize}
\item A {\em logical cut}  occurs when the introduction of the
 connective $\f$ is immediately followed by its
 elimination.  The corresponding reduction rule (denoted by $\be$) is:

$$(\l x M \; N) \tr M[x:=N]$$

\item A {\em classical cut} occurs when $\bot_e$ appears as the left
 premiss of a $\f_e$. The corresponding reduction rule (denoted by
 $\m$) is:

$$(\mu \al M \; N) \tr \mu \al M[\al=_r N]$$

where $M[\al=_r N]$ is obtained by replacing each sub-term of $M$
of the form $(\al \; U)$ by $(\al \; (U \; N))$. This substitution
is called a $\m$-substitution.

\item A {\em symmetric classical cut}  occurs when $\bot_e$ appears as
 the right premiss of a $\f_e$. The corresponding reduction rule (denoted by $\m'$) is:

$$(M \; \mu \al N) \tr \mu \al N[\al=_l M]$$

 where $N[\al=_l M]$
is obtained by replacing each sub-term of $N$ of the form $(\al \;
U)$ by $(\al \; (M \; U))$. This substitution is called a
$\m'$-substitution.
\end{itemize}

\noindent {\bf Remarks}

\begin{enumerate}
  \item It is shown in \cite{Par1} that the $\be\m$-reduction is confluent
but neither $\m\m'$ nor $\be\m'$ is. For example $(\m \al x \, \m
\be y)$ reduces both to $\m \al x$ and to $\m \be y$. Similarly
$(\la z  x \; \m \be y)$ reduces both to $x$ and to $\m \be y$.
  \item  The reductions on terms correspond to the
elimination of cuts on the proofs.
\begin{itemize}
  \item
 The $\be$-reduction is
the usual one.
  \item The $\m$-reduction is as follows. If $M$ corresponds to
a proof of $\perp$ assuming $\al : \neg(A \rightarrow B)$ and $N$
corresponds to a proof of $A$, then $M[\al=_r N]$ corresponds to
the  proof $M$ of $\perp$ assuming $\al : \neg B$ but where, each
time we used the hypothesis $\al:\neg(A \rightarrow B)$ with a
proof $U$ of $A \rightarrow B$ to get $\perp$, we replace this by
the following proof of $\perp$. Use $U$ and $N$ to get a proof of
$B$ and then $\al : \neg B$ to get a proof of $\perp$.
  \item Similarly, the $\m'$-reduction is as follows. If $N$ corresponds
to a proof of $\perp$ assuming $\al : \neg A $ and $M$ corresponds
to a proof of $A \rightarrow B$, then $N[\al=_l M]$ corresponds to
the proof $N$ of $\perp$ assuming $\al : \neg B$ but where, each
time we used the hypothesis $\al:\neg A $ with a proof $U$ of $A$
to get $\perp$, we replace this by the following proof of $\perp$.
Use $U$ and $M$ to get a proof of $B$ and then $\al : \neg B$ to
get a proof of $\perp$.
\end{itemize}

\item Unlike for a $\be$-substitution where, in $M[x:=N]$, the variable $x$ has disappeared it
 is important to note that, in a $\m$ or $\m'$-substitution, the variable $\al$ has not disappeared. Moreover its type has
 changed.  If the type of $N$ is $A$ and, in $ M $, the type of $\al$
 is $\neg (A\rightarrow B)$ it becomes $\neg B$ in $M[\al=_r N]$. If
 the type of $M$ is $A \rightarrow B$ and, in $ N$, the type of $\al$
 is $\neg A$ it becomes $\neg B$ in $N[\al=_l M]$.

\end{enumerate}

\medskip

In the next sections we will study various reductions : the
$\m\m'$-reduction in section \ref{mm'} and the
$\be\m\m'$-reduction in sections \ref{beta}, \ref{s3}. The
following notions will correspond to these reductions.

\begin{definition}\rm
Let $\tr$ be a notion of reduction and $M$ be a term.
\begin{enumerate}
\item   The
 transitive (resp. reflexive and transitive) closure of $\tr$ is denoted by $\tr^+$ (resp. $\tr^*$).

\item If $M$ is in $SN$
  i.e. $M$ has no infinite reduction, $\eta(M)$ will denote the length
of the longest reduction starting from $M$ and $\eta c(M)$ will
denote  $(\eta(M), cxty(M))$.
\item We denote by $N \prec M$ the fact that $N \leq M'$ for some $M'$ such that $M \ras M'$
   and either $M \tr^+ M'$ or $N < M'$. We denote by $\preceq$ the reflexive closure of $\prec$.

\end{enumerate}
\end{definition}

\noindent {\bf Remarks}

- It is easy to check that the relation $\peq$ is transitive and
that $N \preceq M$ iff $N \leq M'$ for some $M'$ such that $M \ras
M'$.

-  If $M \in SN$ and $N \prec M$, then $N \in SN$ and $\eta c(N) <
\eta c(M)$. It follows that the relation $\peq$ is an order on the
set $SN$.

- Many proofs will be done by induction on some $k$-uplet of
integers. In this case the order we consider is the lexicographic
order.

\section{The  $\mu\mu'$-reduction is strongly normalizing}\label{mm'}

In this section we consider the $\m\m'$-reduction, i.e. $M \tr M'$
means $M'$ is obtained from $M$ by one step of the
$\m\m'$-reduction. The main points of the proof of the strong
normalization of $\mu\mu'$ are the following.

\noindent - We first show (cf. lemma \ref{l7}) that a $\m$ or
$\m'$-substitution cannot create a $\m$.

\noindent - It is easy to show (see lemma \ref{lcle}) that if $M
\in SN$ but $M[\sig]\not\in SN$ where $\sig$ is a $\m$ or
$\m'$-substitution, there are an $\al$ in the domain of $\sig$ and
some $M' \prec M$ such that $M'[\sig] \in SN$ and (say $\sig$ is a
$\m$-substitution) $(M'[\sig] \; \sig(\al)) \not \in SN$. This is
sufficient to give a simple proof of the strongly normalization of
the $\mu$-reduction. But this is not enough to do a proof of the
strongly normalization of the $\mu\mu'$-reduction. We need a
stronger (and more difficult) version of this: lemma \ref{crux}
ensure that, if $M[\sig] \in SN$ but $M[\sig][\al=_rP] \not\in SN$
then the real cause of non $SN$ is, in some  sense, $[\al=_rP]$.

\noindent - Having these results, we show, essentially by
induction on $\eta c(M) + \eta c(N)$, that if $M,N \in SN$ then
$(M \; N) \in SN$. The point is that there is, in fact, no deep
interactions between $M$ and $N$ i.e. in a reduct of $(M \; N)$ we
always know what is coming from $M$ and what is coming from $N$.

\begin{definition}\rm
\begin{itemize}
  \item The set of simultaneous substitutions of
  the form $[\al_1=_{s_1}P_1  \ptv$  $ \al_n=_{s_n}P_n]$ where $s_i \in \{l,r\}$ will be denoted by
  $\Sigma$.
  \item  For  $s \in \{l,r\}$, the set of simultaneous substitutions of
  the form $[\al_1=_sP_1 $ ...$\al_n=_sP_n]$  will be denoted by
  $\Sigma_s$.
  \item If $\sig=[\al_1=_{s_1}P_1  \ptv$  $ \al_n=_{s_n}P_n]$, we denote by $dom(\sig)$
  (resp. $Im(\sig)$) the set $\{\al_1, \ptv \al_n\}$
  (resp. $\{P_1, \ptv P_n\}$ ).
\item Let $\sig \in \Sigma$. We say that $\sigma \in SN$ iff for every
$N \in Im(\sig)$, $N \in SN$.
\end{itemize}


\end{definition}

\begin{lemma}\label{ll2}
If $(M \; N) \ras \m \al P$, then either $M \ras \m \al M_1$ and
$M_1[\al =_r N] \ras P$ or $N \ras \m \al N_1$ and $N_1[\al =_l M]
\ras P$.
\end{lemma}

\begin{proof}
By induction  on the length of the reduction $(M \; N) \ras \m \al
P$.
\end{proof}

\begin{lemma}\label{l7}
 Let $M$ be a term and $\sigma \in \Sigma$. If
 $M[\sigma] \ras \m\al P$, then $M\ras \m\al Q$ for some $Q$ such that
 $Q[\sigma] \ras P$.
\end{lemma}

\begin{proof}
By induction on $M$. $M$ cannot be of the form $(\beta \, M')$ or
$\l x \, M'$. If $M$ begins with a $\m$, the result is trivial.
Otherwise $M=(M_1 \; M_2)$ and, by lemma \ref{ll2}, either
$M_1[\sigma] \ras \m\al R$ and $R[\al=_rM_2[\sigma]] \ras P$ or
$M_2[\sigma] \ras \m\al R$ and $R[\al=_lM_1[\sigma]] \ras P$. Look
at the first case (the other one is similar). By the induction
hypothesis $M_1\ras \m\al Q$ for some $Q$ such that $Q[\sigma]
\ras R$ and thus $M \ras \m\al Q[\al=_r M_2]$. Since $Q[\al=_r
M_2][\sigma] = Q[\sigma][\al=_r M_2[\sigma]] \ras R[\al=_r
M_2[\sigma]]\ras P$ we  are done.
\end{proof}

\begin{lemma}\label{l12}
Assume $M,N \in SN$ and $(M \; N) \not\in SN$. Then either $M \ras
\m \al M_1$ and $M_1[\al =_r N] \not\in SN$ or $N \ras \m \be N_1$
and $N_1[\be =_l M] \not\in SN$.
\end{lemma}

\begin{proof}
By induction on $\eta(M)+\eta(N)$. Since $(M \; N) \not\in SN$,
$(M \; N) \tr P$ for some $P$ such that  $P \not \in SN$.  If $P =
(M' \; N)$ where $M \tr M'$ we conclude by the induction
hypothesis since $\eta(M')+\eta(N) < \eta(M)+\eta(N)$. If $P = (M
\; N')$ where $N \tr N'$ the proof is similar. If $M = \mu \al
M_1$  and $P = \mu \al M_1[\al=_rN]$ or  $N = \mu \be N_1$ and $P
= \mu \be N_1[\be=_lM]$
 the result is trivial.
\end{proof}

\begin{lemma}\label{lcle}
Let $M$ be term in $SN$ and $\sig \in \Sigma_s$ be in $SN$. Assume
$M[\sig] \not \in SN$.  Then, for some $(\al \; P) \preceq M$,
$P[\sig] \in SN$ and, if $s=l$ (resp. $s=r$), $(\sig(\al) \,
P[\sig]) \not \in SN$ (resp. $( P[\sig] \, \sig(\al)) \not \in
SN$).
\end{lemma}

\begin{proof}
We only  prove the case $s=l$ (the other one is similar). Let $M_1
\peq M$ be such that $M_1[\sigma] \not\in SN$ and $\eta c(M_1)$ is
minimal.  By the minimality, $M_1$ cannot be $\la x M_2$ or $\m\al
M_2$. It cannot be either $(N_1 \; N_2)$ because otherwise, by the
minimality, the $N_i[\sigma]$ would be in $SN$ and thus, by lemma
\ref{l12} and \ref{l7}, we would have, for example, $N_1 \ras
\m\al N'_1$ and $N'_1[\sigma][\al=_r
N_2[\sigma]]=N'_1[\al=_rN_2][\sigma] \not\in SN$ but this
contradicts the minimality of $M_1$ since $\eta(N'_1[\al=_rN_2]) <
\eta(M_1)$. Then $M_1 = (\al \, P)$ and the the minimality of
$M_1$ implies that $P[\sigma] \in SN$.
\end{proof}

\medskip

\noindent {\bf Remark}

From these results it is easy to prove, by induction on the term,
the strong normalization of the $\mu$-reduction. It is enough to
show that, if $M,N \in SN$, then $(M \; N) \in SN$. Otherwise, we
construct below a sequence $(M_i)$ of terms  and a sequence
$(\sig_i)$ of substitutions  such that, for every $i$,  $\sig_i$
has the form $[\al_1=_rN,...,\al_n=_rN]$, $M_i[\sig_i] \not \in
SN$ and $M_{i+1} \prec M_i \prec M$. The sequence $(M_i)$
contradicts the fact that $M \in SN$. Since $(M \; N) \not \in
SN$,  by lemma \ref{l12}, $M \ras \m \al M_1$ and $M_1[\al =_r N]
\not\in SN$. Assume we have constructed $M_i$ and $\sig_i$. Since
$M_i[\sig_i] \not \in SN$, by lemma \ref{lcle}, there is $M'_i
\prec M_i$ such that $M'_i[\sig_i] \in SN$ and $( M'_i[\sig] \, N)
\not \in SN$. By lemmas \ref{l7} and \ref{l12}, $M'_i \ras \mu \al
M_{i+1}$ and $M_{i+1}[\sig_i+\al=_rN] \not \in SN$.

\medskip

In the remark above, the fact that $(M \; N) \not \in SN$ gives an
infinite $\mu$-reduction in $M$. This not the same for the the
$\mu\mu'$-reduction and, if we try to do the same, the
substitutions we get are more complicated. In particular, it  is
not clear that we get an infinite sequence either of the form $...
\prec  M_{2} \prec M_1 \prec M$ or of the form $... \prec N_{2}
\prec N_1 \prec N$. Lemma \ref{crux} below will give the answer
since it will ensure that, at each step,  we may assume that the
cause of non $SN$ is the last substitution.

\begin{lemma}\label{crux}

Let $M$ be a term and $\sig \in \Sigma_s$. Assume $\delta$ is free in
$M$ but not free in $Im(\sig)$. If $M[\sigma] \in SN$ but
$M[\sigma][\delta=_sP] \not\in SN$, there is $M'\prec M$ and $\sigma'$
such that $M'[\sigma'] \in SN$ and, if $s=r$, $(M'[\sigma'] \;\; P)
\not\in SN$ and, if $s=l$, $(P \;\; M'[\sigma']) \not\in SN$.
\end{lemma}

\begin{proof}
Assume $s=r$ (the other case is similar). Let  $Im(\sigma)= \{N_1,
\ptv N_k\} $. Assume $M,\delta, \sigma, P$ satisfy the hypothesis.
Let ${\cal U}=\{U \; / \; U \peq M\}$ and ${\cal V}=\{V \; / \; V
\peq N_i$ for some $i\}$. Define inductively the sets $\Sigma_m$
and $\Sigma_n$ of substitutions by the following rules:

$\rho \in \Sigma_m$ iff $ \rho = \es$ or $\rho=\rho' + [\be =_r
V[\tau]]$ for some $V \in {\cal V} $, $\tau  \in \Sigma_n$ and
$\rho' \in \Sigma_m$

$\tau \in \Sigma_n$ iff $ \tau = \es$ or $\tau=\tau' + [\al =_l
U[\rho]]$ for some $U \in {\cal U}$, $\rho  \in \Sigma_m$ and
$\tau' \in \Sigma_n$

\noindent Denote by C the conclusion of the lemma, i.e. there is
$M'\prec M$ and $\sigma'$ such that
 $M'[\sigma'] \in SN$, and
 $(M'[\sigma'] \;\; P) \not\in
  SN$.

\noindent We prove something more general.

\noindent (1) Let $U \in {\cal U}$ and $\rho \in \Sigma_m$. Assume
$U[\rho]\in SN$ and $U[\rho][\delta=_rP]\not\in SN$. Then, C
holds.

\noindent (2) Let $V \in {\cal V}$ and $\tau \in \Sigma_n$. Assume
$V[\tau]\in SN$ and $V[\tau][\delta=_rP]\not\in SN$. Then, C
holds.

The conclusion C follows from (1) with $M$ and $\sig$. The properties
(1) and (2) are proved by a simultaneous induction on $\eta
c(U[\rho])$ (for the first case) and $\eta c(V[\tau])$ (for the second
case).

\medskip

Look first at (1)

\noindent - if $U =\l x U'$ or $U=\m\al U'$:  the result follows
 from the induction hypothesis with $U'$ and $\rho$.

\noindent - if $U= (U_1 \; U_2)$: if $U_i[\rho][\delta=_rP]\not\in
SN$ for $i=1$ or $i=2$, the result follows
 from the induction hypothesis with $U_i$ and $\rho$. Otherwise, by lemma \ref{l7} and
 \ref{l12},
 say $U_1 \ras \m\al U'_1$ and, letting $U'=U'_1[\al=_ru_2]$, $U'[\rho][\delta=_rP]\not\in SN$
 and the result follows
 from the induction hypothesis with $U'$ and $\rho$.

\noindent - if $U=(\delta \; U_1)$: if $U_1[\rho][\delta=_rP] \in
SN$, then $M'=U_1$ and $\sig'=\rho[\delta=_rP]$ satisfy the
desired conclusion. Otherwise, the result follows
 from the induction hypothesis with $U_1$ and $\rho$.

\noindent - if $U=(\al \; U_1)$: if $\al \not\in dom(\rho)$ or
$U_1[\rho][\delta=_rP] \not\in SN$, the result follows
 from the induction hypothesis with $U_1$ and $\rho$.  Otherwise,  let
 $\rho(\al)=V[\tau]$. If $V[\tau][\delta=_rP] \not\in SN$, the
 result follows
 from the induction hypothesis with $V$ and $\tau$ (with (2)). Otherwise, by lemma \ref{l7} and
 \ref{l12}, there are two cases to consider.

 - $U_1 \ras \m\al_1 U_2$ and $U_2[\rho'][\delta=_rP] \not\in SN$
 where $\rho' =\rho + [\al_1=_rV[\tau]]$. The result follows from the
 induction hypothesis with $U_2$ and $\rho'$.

 - $V \ras \m\be V_1$ and $V_1[\tau'][\delta=_rP] \not\in SN$ where $\tau'=\tau+[\be=_lU_1[\rho]]$. The result follows
 from the induction hypothesis with $V_1$ and $\tau'$ (with (2)).

\medskip

The case (2) is proved in the same way. Note that, since $\delta$
is not free in the $N_i$, the case $b=(\delta \; V_1)$ does not
appear.
\end{proof}

\begin{theorem}\label{thm1}
Every term is in $SN$.
\end{theorem}

\begin{proof}
By induction on the term. It is enough to show that, if $M,N \in
SN$, then $(M \; N) \in SN$. We prove something more general: let
$\sigma$ (resp. $\tau$) be in $\Sigma_r$ (resp. $\Sigma_l$) and
assume $M[\sigma],N[\tau] \in SN$. Then $(M[\sigma] \; N[\tau])
\in SN$. Assume it is not the case and choose some elements such
that $M[\sigma],N[\tau] \in SN$, $(M[\sigma] \; N[\tau]) \not\in
SN$ and $(\eta(M)+\eta(N), cxty(M)+cxty(N))$ is minimal. By lemma
\ref{l12}, either $M[\sigma] \ras \mu \delta M_1$ and $M_1[\delta
=_r N[\tau]] \not \in SN$ or $N[\tau] \ras \mu \be N_1$ and
$N_1[\be =_l M[\sigma]] \not \in SN$. Look at the first case (the
other one is similar). By lemma \ref{l7}, $M \ras \mu \delta M_2$
for some $M_2$ such that $M_2[\sigma] \ras M_1$. Thus,
$M_2[\sigma][\delta =_r N[\tau]] \not \in SN$. By lemma \ref{crux}
with $M_2, \sig$ and $N[\tau]$, let $M'\prec M_2$ and $\sig'$ be
such that $M'[\sig'] \in SN$, $(M'[\sig'] \; N[\tau]) \not\in SN$.
 This contradicts the
minimality of the chosen elements since $\eta c(M') < \eta c(M)$.
\end{proof}

\section{The simply typed symmetric
$\la\m$-calculus is strongly normalizing}\label{beta}

In this section, we consider the simply typed calculus with the
$\be\m\m'$-reduction i.e. $M \tr M'$ means $M'$ is obtained from
$M$ by one step of the $\be\m\m'$-reduction. To prove the strong
normalization of the $\be\m\m'$-reduction, it is enough to show
that, if $M,N \in SN$, then $M[x:=N]$ also is in $SN$. This is
done by induction on the type of $N$. The proof very much looks
like  the one for the $\m\m'$-reduction and the induction on the
type is used for the cases coming from a $\be$-reduction. The two
new difficulties are the following.

\vspace{.2cm}

\noindent - A $\be$-substitution may create a $\m$, i.e. the fact
that $M[x:=N] \ras \mu \al P$ does not imply that $M \ras \mu \al
Q$. Moreover the $\m$ may come from a  complicated interaction
between $M$ and $N$ and, in particular, the alternation between
$M$ and $N$ can be lost. Let e.g. $M=(M_1 \; (x \; (\la y_1 \la
y_2 \m \al M_4) \; M_2 \; M_3))$ and $N=\la z (z \; N_1)$. Then
$M[x:=N] \tr^* (M_1 \; (\m\al M'_4 \; M_3)) \ras \m\al M'_4[\al
=_rM_3][\al=_lM_1]$. To deal with this situation, we  need to
consider some new kind of $\m\m'$-substitutions (see definition
\ref{def}). Lemma \ref{l8b} gives the different ways in which a
$\m$ may appear. The difficult case in the proof (when a $\m$ is
created and the control between $M$ and $N$ is lost) will be
solved by using a typing argument.

\medskip

\noindent - The crucial lemma (lemma \ref{crux'}) is essentially
the same as the one (lemma \ref{crux}) for the $\m\m'$-reduction
but, in its proof, some cases cannot be proved ``by themselves''
and we need an argument using the types. For this reason, its
proof is done using the additional fact that we already know that,
if $M,N \in SN$ and the type of $N$ is small, then $M[x:=N]$ also
is in $SN$. Since the proof of lemma \ref{thm} is done by
induction on the type, when we will use lemma \ref{crux'}, the
additional hypothesis will be available.

\begin{lemma}\label{lammu}
\begin{enumerate}
\item If $(M \; N) \ras \l x P$, then $M \ras \l y M_1$ and
$M_1[y := N] \ras \l x P$.
\item If $(M \; N) \ras \m \al P$, then either ($M \ras \l y M_1$ and
$M_1[y := N] \ras \m \al P$) or ($M \ras \m \al M_1$ and $M_1[\al
=_r N] \ras P$) or ($N \ras \m \al N_1$ and $N_1[\al =_l M] \ras
P$).
\end{enumerate}
\end{lemma}
\begin{proof}
(1) is trivial. (2) is as in lemma \ref{ll2}.
\end{proof}

\begin{lemma}\label{l8a}
Let $M \in SN$ and $\sig=[x_1:=N_1, ..., x_k:=N_k]$. Assume
$M[\sig] \ras \la y P$. Then, either $M \ras \la y P_1$ and
$P_1[\sig] \ras P$
 or $M \ras  (x_i \; \overrightarrow{Q})$ and $(N_i \;
  \overrightarrow{Q[\sig]})\ras \la y P$.

\end{lemma}
\begin{proof}
By induction on $\eta c(M)$. The only non immediate case is $M=(R
\; S)$. By lemma \ref{lammu}, there is a term $R_1$ such that
$R[\sig] \ras \la z R_1$ and $R_1[z:=S[\sig]] \ras \la y P$. By
the induction hypothesis (since $\eta c(R) < \eta c(M)$), we have
two cases to consider.

(1) $R \ras \la z R_2$ and $R_2[\sig] \ras R_1$, then
$R_2[z:=S][\sig] \ras \la y P$. By the induction hypothesis (since
$\eta(R_2[z:=S])< \eta(M)$),

- either $R_2[z:=S] \ras \la y P_1$ and $P_1[\sig] \ras P$ ; but
then $M \ras \la y P_1$ and we are done.

- or $R_2[z:=S] \ras (x_i \; \overrightarrow{Q})$ and $(N_i \;
\overrightarrow{Q[\sig]})\ras \la y P$, then $M \ras (x_i \;
\overrightarrow{Q})$ and again we are done.

(2) $R \ras (x_i \; \overrightarrow{Q})$ and $(N_i \;
  \overrightarrow{Q[\sig]})\ras \la z R_1$. Then $M \ras (x_i \; \overrightarrow{Q} \;
  S)$ and the result is trivial.
\end{proof}

\begin{definition}\rm\label{def}
\begin{itemize}
  \item An  address is a finite list of symbols in
  $\{l,r\}$. The empty list is denoted by $[]$ and, if $a$ is an
  address and $s \in \{l,r\}$, $[s::a]$ denotes the list obtained
  by putting $s$ at the beginning of $a$.
  \item Let $a$ be an  address and $M$ be a term. The
  sub-term of $M$ at the address $a$ (denoted as $M_a$) is defined
  recursively as follows : if $M=(P \; Q)$ and $a=[r::b]$ (resp. $a=[l::b]$) then
  $M_a=Q_b$ (resp. $P_b$) and undefined otherwise.
  \item Let $M$ be a term and $a$ be an  address
  such that $M_a$ is defined. Then $M\langle a=N \rangle$ is the
  term $M$ where the sub-term $M_a$ has been replaced by
  $N$.
  \item Let $M,N$ be some terms and $a$ be an  address
  such that $M_a$ is defined. Then $N[\al=_aM]$ is the term $N$ in
  which each sub-term of the form $(\al \; U)$ is replaced by
  $(\al \; M\langle a=U\rangle)$.

\end{itemize}
\end{definition}

\noindent {\bf Remarks and examples}

- Let $N = \la x (\al \; \la y (x \; \mu \b (\al \; y)))$, $M = (M_1 \; (M_2
\; M_3))$ and $a=[r::l]$. Then $N[\al=_aM] = \la x (\al \;
(M_1 \; (\la y (x \; \mu \b (\al \;(M_1 \; (y \; M_3)) )) \; M_3)))$.

- Let $M=(P \; ((R \; (x \; T)) \; Q))$ and $a=[r::l::r::l]$. Then
$N[\al=_aM]= N[\al=_r T][\al=_lR][\al=_r Q][\al=_rP]$.

- Note that the sub-terms of a term having an address in the sense
given above are those for which the path to the root consists only
on applications (taking either the left or right son).

 - Note
that $[\al=_{[l]}M]$ is not the same as $[\al=_lM]$ but
$[\al=_lM]$ is the same as $[\al=_{[r]}(M \; N)]$ where $N$ does
not matter. More generally, the term $N[\al=_aM]$ does not depend
of $M_a$.

- Note that  $M\langle a=N \rangle$ can be written as $M'[x_a:=N]$
where $M'$ is the term $M$ in which $M_a$ has been replaced by the
fresh variable $x_a$ and thus (this will be used in the proof of
lemma \ref{thm}) if $M_a$ is a variable $x$, $(\al \; U)[\al=_aM]=
(\al \; M_1[y:=U[\al=_aM]])$ where $M_1$ is the term $M$ in which
the particular occurrence of $x$ at the address $a$ has been
replaced by the fresh name $y$ and the other occurrences of $x$
  remain unchanged.

  \begin{lemma}\label{l30}
Assume $M,N \in SN$ and $(M \; N) \not\in SN$. Then, either ($M
\ras \la y P$ and $P[y:=N] \not\in SN$) or
  ($M \ras \mu \al P$ and $P[\al =_r N] \not\in SN$) or
 ($N \ras \mu \al P$ and $P[\al =_l M] \not\in SN$).
\end{lemma}
\begin{proof}
By induction on $\eta(M)+\eta(N)$.
\end{proof}

\medskip

{\em In the rest of this section, we consider the typed calculus.
To simplify the notations, we do not write explicitly the type
information but, when needed, we denote by $type(M)$ the type of
the term $M$.}

\medskip

\begin{lemma}
If $\G \v M : A$ and $M \ras N$ then $\G \v N : A$.
\end{lemma}

\begin{proof}
Straight forward.
\end{proof}

\begin{lemma}\label{l8b}
Let $n$ be an integer,  $M \in SN$, $\sig=[x_1:=N_1, ...,
x_k:=N_k]$ where $lg(type(N_i))=n$ for each $i$.  Assume $M[\sig]
\ras \mu \al P$. Then,
\begin{enumerate}
  \item either $M \ras \mu \al P_1$ and $P_1[\sig] \ras P$
  \item  or $M \ras  Q$ and, for some $i$, $N_i \ras \m\al N'_i$  and $N'_i[\al=_aQ[\sig]] \ras
P$ for some
 address $a$ in $Q$ such that  $Q_a=x_i$.
  \item or $M \; \ras  Q$,  $Q_a[\sig] \ras \m \al N'$ and
$N'[\al=_aQ[\sig]] \ras P$ for some
 address $a$ in $Q$ such that $lg(type(Q_a)) < n$ .
\end{enumerate}
\end{lemma}
\begin{proof}
By induction on $\eta c(M)$.   The only non immediate case is
$M=(R \; S)$. Since $M[\sig] \ras \mu \al P$, the application
$(R[\sig] \; S[\sig])$ must be reduced. Thus there are three cases
to consider.

\begin{itemize}
  \item It is reduced by a $\m'$-reduction, i.e. there is a term $S_1$
  such that $S[\sig] \ras \mu \al S_1$ and $S_1[\al=_lR[\sig]] \ras
  P$.  By the induction hypothesis: \\ - either $S \ras \m\al Q$ and
  $Q[\sig] \ras S_1$, then $M \ras \m \al Q[\al=_lR]$ and
  $Q[\al=_lR][\sig] \ras P$.\\ - or $S \ras Q$ and, for some $i$, $N_i \ras \m\al N'_i$,
  $Q_a=x_i$ for some address $a$ in $Q$ and $N'_i[\al=_aQ[\sig]] \ras S_1$.
  Then $M \ras (R \; Q)=Q'$ and letting $b=[r::a]$ we have
  $N'_i[\al=_bQ'[\sig]] \ras P$. \\ - or $S \;\ras Q$, $Q_a[\sig] \ras
  \m \al N'$ for some address $a$ in $Q$ such that $lg(type(Q_a)) <
  n$ and $N'[\al=_aQ[\sig]] \ras S_1$. Then $M \ras (R \; Q)=Q'$ and
  letting $b=[r::a]$ we have $N'[\al=_bQ'[\sig]] \ras P$ and
  $lg(type(Q'_b))<n$.
  \item It is reduced by a $\m$-reduction. This
  case is similar to the previous one.
   \item It is reduced by a
  $\be$-reduction, i.e. there is a term $U$ such that $R[\sig] \ras
  \la y U$ and $U[y:=S[\sig]] \ras \mu \al P$. By lemma \ref{l8a},
  there are two cases to consider. \\ - either $R \ras \la y R_1$ and
  $R_1[\sig][y:=S[\sig]]=R_1[y:=S][\sig] \ras \mu \al P$. The result
  follows from the induction hypothesis sine $\eta(R_1[y:=S]) <
  \eta(M)$.\\ - or $R \ras (x_i \; \overrightarrow{R_1})$. Then $Q=(x_i \;
\overrightarrow{R_1} \; S)$ and $a=[]$ satisfy the desired
conclusion since then $lg(type(M)) <n$. \qed
\end{itemize}
\end{proof}

\begin{definition}\rm
Let $A$ be a type. We denote by $\Sigma_A$ the set of
substitutions of the form $[\al_1=_{a_1}M_1, ...,
  \al_n=_{a_n}M_n]$  where the type of the $\al_i$ is $\neg A$.
\end{definition}

\noindent {\bf Remark}

 Since in such  substitutions  the type of
the variables changes, when  we consider the term $N[\sig]$ where
$\sig \in \Sigma_A$, we mean that the type  of the $\al_i$ is $A$
in $N$ i.e. before the substitution. Also note that considering
$N[\al=_aM]$ implies that the type of $M_a$ is $A$.

\begin{lemma}\label{crux'}
Let $n$ be an integer and $A$ be a type such that $lg(A)=n$. Let
$N,P$ be terms and $\tau \in \Sigma_A$. Assume that,
\begin{itemize}
  \item  for every $M,N \in SN$ such that $lg(type(N))<n$,
$M[x:=N] \in SN$.
  \item $N[\tau] \in SN$ but $N[\tau][\delta =_aP] \not \in SN$.
  \item $\delta$ is free and has type $\neg A$ in $N$ but $\delta$ is not free in
$Im(\tau)$.
\end{itemize}

\noindent Then,  there is $N'\prec N$ and $\tau' \in \Sigma_A$
such that
 $N'[\tau'] \in SN$ and
 $P\langle a = N'[\tau']\rangle\not\in
  SN$.

\end{lemma}
\begin{proof}
Essentially as in lemma \ref{crux}. Denote by (H) the first
assumption i.e. for every $M,N \in SN$ such that $lg(type(N))<n$,
$M[x:=N] \in SN$.

Let $\tau =[\al_1=_{a_1}M_1, ...,
  \al_n=_{a_n}M_n]$, ${\cal U}=\{U \; / \; U \peq N\}$ and
${\cal V}=\{V \; / \; V \peq M_i$ for some $i\}$. Define
inductively the sets $\Sigma_m$ and $\Sigma_n$ of substitutions by
the following rules:

$\rho \in \Sigma_n$ iff $ \rho = \es$ or $\rho=\rho'  + [\al =_a
V[\sig]]$ for some $V \in {\cal V} $, $\sig  \in \Sigma_m$, $\rho'
\in \Sigma_n$ and $\al$ has type $\neg A$.

$\sig \in \Sigma_m$ iff $ \sig = \es$ or $\sig=\sig' + [x:=
U[\rho]]$ for some $U \in {\cal U}$, $\rho  \in \Sigma_n$, $\sig'
\in \Sigma_m$ and $x$ has type $A$.

\noindent Denote by C the conclusion of the lemma. We prove
something more general.

\noindent (1) Let $U \in {\cal U}$ and $\rho \in \Sigma_n$. Assume
$U[\rho]\in SN$ and $U[\rho][\delta=_aP]\not\in SN$. Then, C
holds.

\noindent (2) Let $V \in {\cal V}$ and $\sig \in \Sigma_m$. Assume
$V[\sig]\in SN$ and $V[\sig][\delta=_aP]\not\in SN$. Then, C
holds.

The conclusion C follows from (1) with $N$ and $\tau$. The properties
(1) and (2) are proved by a simultaneous induction on $\eta
c(U[\rho])$ (for the first case) and $\eta c(V[\tau])$ (for the second
case).

 The proof is as in
lemma \ref{crux}. The new case to consider is, for $V[\sig]$, when
$V=(V_1 \; V_2)$ and $V_i[\sig][\delta=_aP]\in SN$.

\medskip

\noindent - Assume first the interaction between $V_1$ and $V_2$
is a $\be$-reduction. If $V_1 \ras \l x V'_1$, the result follows
from the induction hypothesis with $V'_1[x:=V_2][\sig]$.
Otherwise, by lemma \ref{l8a}, $V_1 \ras (x \;
\overrightarrow{W})$. Let $\sig(x)=U[\rho]$. Then  $(U[\rho]
\;\overrightarrow{W}[\sig]) \ras \l y Q$ and
$Q[y:=V_2[\sig]][\delta=_aP]\not\in SN$. But, since the type of $x$ is
$A$, the type of $y$ is less than $A$ and since $Q[\delta=_aP]$ and
$V_2[\sig][\delta=_aP]$ are in $SN$ this contradicts (H).

\noindent - Assume next the interaction between $V_1$ and $V_2$ is
a $\m$ or $\m'$-reduction.  We consider only the case $\m$ (the
other one is similar). If $V_1 \ras \m \al V'_1$, the result
follows from the induction hypothesis with
$V'_1[\al=_rV_2][\sig]$. Otherwise, by lemma \ref{l8b}, there are
two cases to consider.

- $V_1 \ras Q$, $Q_c=x$ for some address $c$ in $Q$ and $x \in
dom(\sig)$, $\sig(x)= U[\rho]$, $U \ras \m\al U_1$ and
$U_1[\rho][\al=_cQ[\sig]] [\al=_r V_2[\sig]][\delta=_aP]\not\in
SN$. Let $V'=(Q \; V_2)$ and $b=l::c$. The result follows then from
the induction hypothesis with $U_1[\rho']$ where $\rho'=\rho +
[\al=_bV'[\sig]]$.

- $V_1 \; \ras Q$, $Q_c[\sig][\delta=_aP] \ras \m \al R$ for some
address $c$ in $Q$ such that $lg(type(Q_c)) < n$,
$R[\al=_cQ[\sig][\delta=_aP]] [\al=_r
V_2[\sig][\delta=_aP]]\not\in SN$. Let $V'=(Q' \; V_2)$ where $Q'$
is the same as $Q$ but $Q_c$ has been replaced by a fresh variable
$y$ and $b=l::c$. Then $R[\al=_bV'[\sig][\delta=_aP]] \not \in
SN$. Let $R'$ be such that $R' \prec R$, $R'
[\al=_bV'[\sig][\delta=_aP]] \not \in SN$ and $\eta c(R')$ is
minimal. It is easy to check that $R'=(\al \; R'')$,
$R''[\al=_bV'[\sig][\delta=_aP]] \in SN$ and
$V'[\sig'][\delta=_aP] \not\in SN$ where $\sig' =\sig
+y:=R''[\al=_bV'[\sig]]$. If $V'[\sig][\delta=_aP] \not\in SN$, we
get the result by the induction hypothesis since $\eta c(V'[\sig])
< \eta c (V[\sig])$. Otherwise this contradicts the assumption (H)
since $V'[\sig][\delta=_aP], R''[\al=_bV'[\sig][\delta=_aP]] \in
SN$, $V'[\sig][\delta=_aP][y:= R''[\al=_bV'[\sig][\delta=_aP]]]
\not\in SN$ and the type of $y$ is less than $n$.

\end{proof}

\begin{lemma}\label{thm}
If $M,N \in SN$, then $M[x:=N] \in SN$.
\end{lemma}
\begin{proof}
We prove something a bit more general: let $A$ be a type, $M, N_1,
..., N_k$ be terms and $\tau_1, ..., \tau_k $ be substitutions in
$\Sigma_A$. Assume that, for each $i$,  $N_i$ has type $A$ and
$N_i[\tau_i] \in SN$. Then $M[x_1:=N_1[\tau_1], \ptv
x_k:=N_k[\tau_k]] \in SN$. This is proved by induction on $(lg(A),
\eta(M), cxty(M)$, $\Sigma \; \eta(N_i), \Sigma \; cxty(N_i))$
where, in $\Sigma \; \eta(N_i)$ and $\Sigma \; cxty(N_i)$, we
count each occurrence of the substituted variable. For example if
$k=1$ and $x_1$ has $n$ occurrences, $\Sigma \; \eta(N_i)=n.
\eta(N_1)$.

 If $M$ is $\lambda y M_1$ or $(\alpha \; M_1)$ or $\mu\alpha M_1$
or a variable, the result is trivial.
  Assume then that $M=(M_1 \;
M_2)$. Let $\sigma = [x_1:=N_1[\tau_1], \ptv x_k:=N_k[\tau_k]]$.
By the induction hypothesis, $M_1[\sigma], M_2[\sigma] \in SN$. By
lemma \ref{l30} there are 3 cases to consider.

\begin{itemize}
  \item  $M_1[\sigma] \ras \lambda y P$ and $P[y:=M_2[\sigma]] \not\in
SN$. By lemma \ref{l8a}, there are two cases to consider.

\begin{itemize}
  \item $M_1 \ras \lambda y Q$ and $Q[\sigma] \ras P$. Then
  $Q[y:=M_2][\sigma]= Q[\sigma][y:=M_2[\sigma]] \ras
  P[y:=M_2[\sigma]]$ and, since $\eta(Q[y:=M_2]) < \eta(M)$, this contradicts the induction hypothesis.
  \item $M_1 \ras (x_i \; \overrightarrow{Q})$ and $(N_i \;
  \overrightarrow{Q[\sigma]}) \ras \lambda y P$. Then, since the type of $N_i$ is $A$,
 $lg(type(y))<lg(A)$. But $P, M_2[\sigma]
  \in SN$ and $P[y:=M_2[\sigma]] \not\in
SN$. This contradicts the induction hypothesis.

\end{itemize}

  \item $M_1[\sigma] \ras \mu\alpha P$ and $P[\alpha=_rM_2[\sigma]]
  \not\in SN$. By lemma \ref{l8b}, there are three cases to consider.

\begin{itemize}
  \item $M_1 \ras \mu\alpha Q$ and $Q[\sigma] \ras P$. Then, $Q[\alpha=_rM_2][\sigma]= Q[\sigma][\alpha=_rM_2[\sigma]] \ras
  P[\alpha=_rM_2[\sigma]]$ and, since $\eta(Q[\alpha=_rM_2]) < \eta(M)$, this contradicts the induction hypothesis.
  \item $M_1 \ras  Q$, $N_i [\tau_i] \ras \m\alpha L'$ and  $Q_a=x_i$ for some
 address $a$ in $Q$ such that $L'[\alpha=_aQ[\sigma]] \ras
P$ and thus $L'[\al=_bM'[\sigma]]  \not\in SN$ where $b=(l::a)$ and
$M'=(Q \; M_2)$.

By lemma \ref{l7}, $N_i \ras \mu\al L$ and $L[\tau_i] \ras L'$.
Thus, $L[\tau_i][\al=_bM'[\sigma]]  \not\in SN$. By lemma
\ref{crux'}, there is $L_1 \prec L$ and $\tau'$ such that
$L_1[\tau'] \in SN$ and $M'[\sigma]\langle
b=L_1[\tau']\rangle\not\in SN$. Let $M''$ be $M'$ where the
variable $x_i$ at the address $b$ has been replaced by the fresh
variable $y$ and let $\sig_1=\sig + y:=L_1[\tau']$. Then
$M''[\sig_1]=M'[\sigma]\langle b=L_1[\tau']\rangle\not\in SN$.

If $M_1 \ra^+  Q$ we get a contradiction from the induction
hypothesis since $\eta(M'') < \eta(M)$. Otherwise, $M''$ is the
same as $M$ up to the change of name of a variable and $\sigma_1$
differs from $\sigma$ only at the address $b$. At this address,
$x_i$ was substituted in $\sigma$ by $N_i[\tau_i)$ and in
$\sigma_1$ by $L_1[\tau']$ but $\eta c(L_1) < \eta c(N_i)$ and
thus we get a contradiction from the induction hypothesis.

  \item $M \ras  Q$,  $Q_a[\sigma] \ras \m \al L$ for some
 address $a$ in $Q$ such that $lg(type(Q_a)) < lg(A)$  and
$L[\al=_aQ[\sigma]] \ras P$. Then, $L[\al=_bM'[\sigma]] \not \in
SN$ where $b=[l::a]$ and $M'=(Q \; M_2)$.

By lemma \ref{crux'}, there is an $ L'$ and  $\tau'$ such that
$L'[\tau'] \in SN$ and $M'[\sigma]\langle b=L'[\tau'] \rangle\not
\in SN$. Let $M''$ be $M'$ where the variable $x_i$ at the address
$b$ has been replaced by the fresh variable $y$. Then $M''[\sig][
y:=L'[\tau']]=M'[\sigma]\langle b=L'[\tau']\rangle\not\in SN$.

But $\eta(M'') \leq \eta(M)$ and $cxty(M'') < cxty(M)$ since,
because of its type, $Q_a$ cannot be a variable and thus, by the
induction hypothesis, $M''[\sigma]\in SN$. Since $M''[\sigma ][
y:=L'[\tau']]\not \in SN$ and $lg(type(L')) < lg(A)$, this
contradicts the induction hypothesis.

\end{itemize}

  \item $M_2[\sigma] \ras \mu\alpha P$ and $P[\alpha=_lM_1[\sigma]]
  \not\in SN$. This case is similar to the previous one. \qed
\end{itemize}
\end{proof}

\begin{theorem}
Every typed term is in $SN$.
\end{theorem}
\begin{proof}
 By induction on the term. It is enough to show that if $M,N \in SN$, then $(M \; N)
\in SN$. Since $(M \; N)=(x \; y)[x:=M][y:=N]$ where $x,y$ are
fresh variables, the result follows by applying  theorem \ref{thm}
twice and the induction hypothesis.
\end{proof}

\section{Why the usual candidates do not work ?}\label{s3}

In \cite{Par2}, the proof of the strong normalization of the
$\la\m$-calculus is done by using the {\em usual} (i.e. defined
without a fix-point operation) candidates of reducibility. This
proof could be easily extended to the symmetric $\la\m$-calculus
if we knew the following properties for the un-typed calculus:

\begin{enumerate}
  \item If $N$ and $(M[x:=N] \; \overrightarrow{P})$ are in $SN$,
  then so is $(\la x M \; N \; \overrightarrow{P})$.
  \item If $N$ and $(M[\al=_rN] \; \overrightarrow{P})$ are in $SN$,
  then so is $(\m \al M \; N \; \overrightarrow{P})$.
  \item If $\overrightarrow{P}$ are in $SN$, then so is $(x \; \overrightarrow{P})$.
\end{enumerate}

These properties are easy to show for the $\be\m$-reduction but
they were not known for the $\be\m\m'$-reduction.

The properties (1) and (2) are false. Here is a counter-example.
Let $M_0=\la x (x \; P \; \underline{0})$ and $M_1=\la x (x \; P
\;
 \underline{1})$ where $\underline{0} = \la x \la y y$, $\underline{1} = \la x \la
 y x$,  $\Delta = \la x
(x \; x)$ and $P=\la x \la y \la z \;(y \; (z \;\underline{1} \;
\underline{0}) \; (z \;
 \underline{0} \; \underline{1}) \;
\la d \underline{1}   \; \Delta  \;  \Delta)$. Let $M = \la f (f
\; (x \; M_1) \; (x \; M_0))$, $M' = \la f (f \; (\be \; \la x (x
\; M_1)) \; (\be \; \la x (x \; M_0)))$
 and $N= (\al \; \la z(\al \;z))$. Then,
 \begin{itemize}
  \item $M[x:=\m\al N] \in SN$ but $(\la x M \; \m\al N) \not\in SN$.
  \item $M'[\be=_r\m\al  N] \in SN$ but $(\m \be  M' \; \m\al  N) \not\in SN$.
\end{itemize}

This comes from the fact that $(M_0 \; M_0)$ and $(M_1 \; M_1)$
are in $SN$ but $(M_1 \; M_0)$ and $(M_0 \; M_1)$ are not in $SN$.
More details can be found in \cite{craco}.

\medskip

The third property is true and its  proof is essentially the same
as the one of the strong normalization of $\mu\mu'$.   This comes
from the fact that, since $(x \; M_1 ... M_n)$ never reduces to a
$\la$, there is no ``dangerous'' $\be$-reduction. In particular,
the $\be$-reductions we have to consider in the proofs of the
crucial lemmas, are uniquely those that appear in the reductions
$M\peq M'$. We give this proof below.

\begin{lemma}\label{lambda}
The term $(x \; M_1 \pts  M_n)$ never reduces to a term of the
form $\la y M$.
\end{lemma}
\begin{proof}
By induction on $n$. Use lemma \ref{lammu}.
\end{proof}

\begin{definition}\rm
\begin{itemize}
  \item Let $M_1, ..., M_n$ be terms and $1\leq i \leq n$. Then, the term
$M$ in which every sub-term of the form $(\al \; U)$ is replaced
by $(\al \; (x \; M_1 \pts M_{i-1} \; U \; M_{i+1}$ ... $ M_n))$
will be denoted by  $M[\al=_i (M_1 \pts M_n)]$.
  \item We will  denote by $\Sigma_x$ the set of simultaneous substitutions of the form
$[\al_1=_{i_1}(M^1_1 \pts M^1_n),...,\al_k=_{i_k}(M^k_1 \pts
  M^k_n)]$ .

\end{itemize}
\end{definition}

\noindent {\bf Remark}

These substitutions are special cases of the one defined in
section \ref{beta} (see definition \ref{def}). For example
$M[\al=_2 (M_1 \; M_2 \; M_3)]=M[\al=_l (x \; M_1)][\al=_r
M_3]=M[\al=_a (x \; M_1 \; M_2 \; M_3)] $ where $a=[l::r]$.

\begin{lemma}\label{l21}
Assume $(x \; M_1 \pts  M_n) \ras \m\al M$. Then, there is an $i$
such that $M_i \ras \m\al P$ and $P[\al=_i (M_1 \pts M_n)] \ras
M$.
\end{lemma}

\begin{proof}
By induction on $n$. \\
- $n=1$. By lemma \ref{lammu}, $M_1 \ras \m\al P$ and
$P[\al=_l x] = P[\al=_1 (M_1)] \ras M$. \\
- $n \geq 2$. Assume $(x \; M_1 \pts  M_{n-1} \,M_n) \ras \m\al
M$. By lemmas \ref{lammu} and \ref{lambda},

- either $(x \; M_1 \pts  M_{n-1}) \ras \m\al N$ and $N[\al=_rM_n]
  \ras M$. By the induction hypothesis,  there is an $i$
such that $M_i \ras \m\al P$ and $P[\al=_i (M_1 \pts M_{n-1})]
\ras N$. Then  $P[\al=_i (M_1 \pts M_{n-1}\,M_n)] = P[\al=_i (M_1
\pts
  M_{n-1})][\al=_rM_n] \ras N[\al=_rM_n] \ras M$.

- or $M_n \ras \mu \al N$ and $N[\al=_l(x \; M_1 \pts  M_{n-1})]
\ras
  M$. Then $N[\al=_l(x \; M_1 \pts  M_{n-1})]$ = $N[\al=_n(M_1 \pts
  M_{n-1}\,M_n)] \ras M$.
\end{proof}

\begin{lemma}\label{l22}
Assume $M_1,...,M_n \in SN$ and $(x \; M_1 \pts M_n)\not\in SN$.
Then,  there is an $1\leq i \leq n$ such that $M_i \ras \mu \al
\;U$ and $U[\al=_i (M_1 \pts M_n)] \not\in SN$.
\end{lemma}
\begin{proof}
Let $k$ be the least such that $(x \; M_1 \pts M_{k-1}) \in SN$
and $(x \; M_1 \pts M_k) \not\in SN$. By lemmas \ref{l30} and
\ref{lambda},

- either $M_k \ras \m\al U$ and $U[\al=_l (x \; M_1 \pts M_{k-1})]
\not\in SN$. Then, $i=k$ satisfies the desired property since
$U[\al=_k (M_1 \pts M_n)]=U[\al=_l (x \; M_1 \pts M_{k-1})][\al=_r
  M_{k+1}] ...[\al=_r M_n]$.

- or $(x \; M_1 \pts M_{k-1}) \ras \m\al P$ and $P[\al=_r M_k]
\not\in SN$. By lemma \ref{l21}, let $i \leq k-1$ be such that
that $M_i \ras \m\al U$ and $U[\al=_i (M_1 \pts M_{k-1})] \ras P$.
Then $U[\al=_i (M_1 \pts M_n)] \not\in SN$ since $U[\al=_i (M_1
\pts
  M_n)]= U[\al=_i (M_1 \pts M_{k-1})][\al=_r M_k][\al=_r M_{k+1}]
...[\al=_r M_n]$ reduces to $P[\al=_r M_k][\al=_r M_{k+1}]
...[\al=_r M_n]$.
\end{proof}

\begin{lemma}\label{l7a}
 Let $M$ be a term and $\sigma \in \Sigma_x$. If $M[\sigma] \ras \m\al
 P$ (resp. $M[\sigma] \ras \l x P$) , then $M\ras \m\al Q$
 (resp. $M\ras \l x Q$) for some $Q$ such that $Q[\sigma] \ras P$.
\end{lemma}

\begin{proof}
As in lemma \ref{l7}.
\end{proof}

\begin{lemma}\label{crux''}

Let $M$ be a term and $\sigma \in \Sigma_x $. Assume $\delta$ is
free in $M$ but not free in $Im(\sig)$. If $M[\sigma] \in SN$ but
$M[\sigma][\delta=_i(P_1 ... P_n)] \not\in
  SN$, there is $M'\prec M$ and $\sigma'$ such that
 $M'[\sigma'] \in SN$ and $(x \;
 P_1 ... P_{i-1} \; M'[\sigma']\; P_{i+1} ...P_n)\not\in
  SN$.
\end{lemma}

\begin{proof}
 As in lemma \ref{crux}.
\end{proof}

\begin{theorem}\label{xt1}
Assume $M_1,..., M_n$ are in $SN$. Then $(x \; M_1 \pts M_n) \in
SN$.
\end{theorem}
\begin{proof}
We prove a more general result: Let  $M_1,..., M_n$ be terms and
$\sig_1,..., \sig_n$ be in $\Sigma_x$.  If $M_1[\sig_1], \ptv
M_n[\sig_n] \in SN$, then $(x \; M_1[\sig_1] \pts M_n[\sig_n]) \in
SN$. The proof is done exactly as in theorem \ref{thm1} using
lemmas \ref{l22}, \ref{l7a} and \ref{crux''}.
\end{proof}

\section{Future work}\label{s7}

\begin{itemize}
  \item Parigot has introduced other simplification rules in the
$\l\m$-calculus. They are as follows : $ (\al \; \m\be M)
\rightarrow_{\rho} M[\be:=\al]$ and, if $\al$ is not free in $M$,
$ \m\al (\al \; M) \rightarrow_{\theta} M$. It would be
interesting to extend our proofs to these reductions. The rule
$\theta$ causes no problem since it is strongly normalizing and it
is easy to see that this rule can be postponed (i.e. if $M
\rightarrow^*_{\be\m\m'\rho\theta} M_1$ then $M
\rightarrow^*_{\be\m\m'\rho} M_2 \rightarrow^*_{\theta} M_1$ for
some $M_2$). However it is not the same for the rule $\rho$ which
cannot be postponed. Moreover a basic property (if $M[\al=_sN]
\ras \m\be P$, then $M \ras \m \be Q$ for some $Q$ such that
$Q[\al=_s N] \ras P$) used in the proofs is no more true if the
$\rho$-rule is used. It seems that, in this case,  the $\m$ can
only come either from $M$ or from $N$ i.e. without deep
interaction between $M$ and $N$ and thus that our proofs can be
extended to this case but, due to the lack of time, we have not
been able to check the details.
  \item We believe that our technique, will allow to give explicit bounds for
  the length of the reductions of a typed term. This is a goal we
  will try to manage.

\end{itemize}

\end{document}